\documentclass{commat}

\title{%
Quasi Yamabe Solitons on 3-Dimensional Contact Metric Manifolds with $Q\varphi=\varphi Q$
    }
\author{%
V. Venkatesha and H. Aruna Kumara
    }

\affiliation{
    \address{V. V. --
    Department of Mathematics.  Kuvempu University, Shankaraghatta.  Shivamogga, Karnataka-577 451.  India.}
    \email{vensmath@gmail.com}
    \address{H. A. K. --
    Department of Mathematics.  BMS Institute of Technology and Management, Bangalore-64, Karnataka, India.}
    \email{arunakumara@bmsit.in}
    }

\abstract{%
In this paper we initiate the study of quasi Yamabe soliton on 3-dimensional contact metric manifold with $Q\varphi=\varphi Q$ and prove that if a 3-dimensional contact metric manifold $M$ such that $Q\varphi=\varphi Q$ admits a quasi Yamabe soliton with non-zero soliton vector field $V$ being point-wise collinear with the Reeb vector field $\xi$, then $V$ is a constant multiple of $\xi$, the scalar curvature is constant and the manifold is Sasakian. Moreover, $V$ is Killing. Finally, we prove that if $M$ is a 3-dimensional compact contact metric manifold such that $Q\varphi=\varphi Q$ endowed with a quasi Yamabe soliton, then either $M$ is flat or soliton is trivial.
    }


\VOLUME{30}
\NUMBER{1}
\firstpage{191}
\DOI{https://doi.org/10.46298/cm.9695}

\begin{paper}

\section{Introduction }
The concept of Yamabe flow was introduced by Hamilton at the same time as the Ricci flow \cite{Hamilton} as a tool to construct Yamabe metrics on compact Riemannian manifold. On a smooth Riemannian manifold, Yamabe flow can be defined as the evolution of the Riemannian metric $g_0$ in time $t$ to $g=g(t)$ by means of the equation
\begin{align*}
	\frac{\partial}{\partial t}g=-rg, \quad g(0)=g_0,
\end{align*} 
where $r$ denotes the scalar curvature corresponds to $g$. The significance of Yamabe flow lies in the fact that it is a natural geometric deformation to metrics of constant scalar curvature. The Yamabe flow is equivalent to Ricci flow when the dimension is 2. However in dimension greater than 2, the Yamabe flow and Ricci flow are different in their properties, since the first one preserves the conformal class of the metric but the Ricci flow does not in general. A Yamabe solitons corresponds to self-simillar solution of the Yamabe flow and it is defined on Riemannian manifold $(M,g)$ by a vector field $V$ satisfying the equation \cite{Barbosa}
\begin{align}
	\frac{1}{2}\mathcal{L}_Vg=(r-\lambda)g,
\end{align}
where $L_V$ denotes the Lie derivative in the direction of the vector field $V$ and $\lambda$ is a real number. A Yamabe soliton is said to be expanding, steady,
or shrinking, respectively, if $\lambda<0$, $\lambda=0$, or $\lambda>0$. In \cite{Chen}, Chen and Deshmukh introduce the concept of quasi Yamabe soliton. A smooth Riemannian manifold $(M,g)$ is said to be a quasi Yamabe soliton if it admits a vector field $V$ such that 
\begin{align}
	\label{1.2}\frac{1}{2}\mathcal{L}_Vg=(r-\lambda)g+\frac{1}{m}V^\#\otimes V^\#,
\end{align} 
for some constant $m$, where $V^\#$ is the dual 1-form of $V$. The vector field $V$ is called a soliton vector field for quasi Yamabe soliton. Note that if $m\rightarrow\infty$, the equation \eqref{1.2} defines a Yamabe soliton. So, the quasi Yamabe solitons are natural extension of the Yamabe solitons. When the vector field $V$ is gradient of smooth function $f:M\rightarrow\mathbb{R}$, the equation \eqref{1.2} becomes
\begin{align}
	\label{1.3}Hessf=(r-\lambda)g+\frac{1}{m}df\otimes df,
\end{align}
where $Hessf$ denotes the Hessian of $f$ and in this case $f$ is called potential function of the quasi Yamabe soliton and $g$ is called quasi Yamabe gradient soliton and this concept was introduced by Huang and Li \cite{Huang}. When $f$ is constant function, we say that $(M,g)$ is a trivial quasi Yamabe gradient soliton.
\par There are several works on Yamabe flow and quasi Yamabe flow by Barbosa and Ribeiro \cite{Barbosa}, Brendle \cite{Brendle1}, Bruchard et al. \cite{Bruchard}, Hsu \cite{Hsu1}, Daskalopoulos and Sesum \cite{Daskalopoulos}, Neto and Pina \cite{Neto}, Ma and Cheng \cite{Ma}, Ma and Miquel \cite{Ma1}, Cao et al. \cite{Cao}, Blaga \cite{Blaga}, Chen and Deshmukh \cite{Chen}, Pirhadi and Razavi \cite{Pirhadi}, Wang \cite{Wang}, Yang and Zhang \cite{Yang} and many others. In the background of contact geometry, Yamabe solitons were first studied by Sharma \cite{Sharma} in 3-dimensional Sasakian manifold. Then these are studied by Wang \cite{YWang}, Suh and Mandal \cite{Suh} and Erken \cite{Erken}. 

Recently, the present first author with Naik studied the Yamabe soliton on 3-dimen\-sion\-al contact metric manifold with $Q\varphi=\varphi Q$ in \cite{Venkatesha} and generalized the result of Sharma \cite{Sharma}. In Sasakian manifold the Ricci operator $Q$ commutes with the fundamental collineation $\varphi$; but in contact metric manifold this commutativity need not hold. In this paper, we start a study of quasi Yamabe soliton on 3-dimensional contact metric manifold with $Q\varphi=\varphi Q$ and first we prove the following:
\begin{theorem}\label{t1.1}
	If a 3-dimensional contact metric manifold $M$ such that $Q\varphi=\varphi Q$ admits a quasi Yamabe soliton with non-zero soliton vector field $V$ being point-wise collinear with the Reeb vector field $\xi$, then $V$ is a constant multiple of $\xi$, the scalar curvature is constant and the manifold is Sasakian. Moreover, the soliton vector field is Killing. 
\end{theorem}
Next, we consider a 3-dimensional compact contact metric manifold with $Q\varphi=\varphi Q$ and prove
\begin{theorem}\label{t1.2}
	If a 3-dimensional compact contact metric manifold $M$ with $Q\varphi=\varphi Q$ admits a quasi Yamabe soliton, then $M$ is a flat or soliton is trivial.
\end{theorem}

\section{Preliminaries}
In the present section, we shall collect some fundamental results on contat metric manifold (for more details see \cite{Blair}) which will be used to prove our results. A detailed study on contact metric manifolds under various condition are available in \cite{Naik}, \cite{Ghosh}, \cite{Kumara}, \cite{Venkatesha}, \cite{YWang} and references therein.
\par Let $M$ be a smooth differentiable manifold of dimension $(2n+1)$. If on $M$ there exist a (1,1)-type tensor field $\varphi$, vector field $\xi$ and 1-form $\eta$ such that 
\begin{align}
	\label{2.1}\varphi^2=-I+\eta\otimes\xi, \qquad \eta(\xi)=1,
\end{align}
then we say that the triplet $(\varphi,\xi,\eta)$ is an almost contact structure on $M$, and $\xi$ is called the characteristic or the Reeb vector field. It follows from \eqref{2.1} that $\varphi \xi=0$, $\eta\circ\varphi=0$ and $rank(\varphi)=2n$. In general, a smooth manifold $M$ with such a structure is called an almost contact manifold which is denoted by $(M,\varphi,\xi,\eta)$. It is well known that a smooth manifold $M$ admits an almost contact structure if and only if the structure group of the tangent bundle of $M$ reduces to $U(n)\times 1$. If an almost contact manifold admits a Riemannian metric $g$ such that
\begin{align}
	g(\varphi X,\varphi Y)=g(X,Y)-\eta(X)\eta(Y),
\end{align}
for all vector fields $X,Y$ on $M$, then $g$ is said to be compatible with the almost contact structure and $M$ together with $(\varphi,\xi,\eta,g)$ is called an almost contact metric manifold. An almost contact metric manifold $(M,\varphi,\xi,\eta,g)$ is said to be contact metric manifold if $d\eta(X,Y)=g(X,\varphi Y)$, for all vector fields $X,Y$ on $M$. We recall two self-adjoint operators $h=\frac{1}{2}\mathcal{L}_\xi \varphi$, where $\mathcal{L}$ is the usual Lie derivative, and $l=R(\cdot,\xi)\xi$, where $R$ is Riemannian curvature. The tensor $h$ and $l$ satisfy $h\varphi=-\varphi h$, $h\xi=0$, $trh=tr\varphi h=0$ and $l\xi=0$. Furthermore, we also have
\begin{align}
	\label{2.3}\nabla_X \xi=-\varphi X-\varphi hX,\\
	\label{2.4}trl=Ric(\xi,\xi)=2n-trh^2,
\end{align}
where $\nabla$ is the Levi-Civita connection of $g$ and $Ric$ is the Ricci tensor, is defined by $Ric(X,Y)=g(QX,Y)$.

A contact metric manifold is said to be $K$-contact if $\xi$ is Killing with respect to $g$, equivalently, $h=0$, or $trl=2n$. An almost contact metric structure on $M$ is said to be normal if the almost complex structure $J$ on the product manifold $M\times\mathbb{R}$ defined by
\begin{align*}
	J\left(X,f\frac{d}{dt}\right)=\left(\varphi X-f\xi,\eta(X)\frac{d}{dt}\right),
\end{align*}
where $X$ denotes the vector field tangent to $M$, $t$ is the coordinate of $\mathbb{R}$ and $f$ is a smooth function defined on the product manifold $M\times\mathbb{R}$, is integrable. According to Blair \cite{Blair}, the normality of an almost contact structure is expressed by $[\varphi,\varphi]=-2d\eta\otimes\xi$, where $[\varphi,\varphi]$ denotes the Nijenhuis tensor of $\varphi$ defined by
\begin{align*}
	[\varphi,\varphi](X,Y)=\varphi^2[X,Y]+[\varphi X, \varphi Y]-\varphi[\varphi X,Y]-\varphi[X, \varphi Y],
\end{align*}
for any vector fields $X,Y$ on $M$. A normal contact metric manifold
is called a Sasakian manifold. Sasakian manifolds are $K$-contact and 3-dimensional
$K$-contact manifolds are Sasakian.

On 3-dimensional contact metric manifold with $Q\varphi=\varphi Q$, the following relations hold (see \cite{BKS}):
\begin{align}
	\nonumber R(X,Y)Z=&(\frac{r}{2}-trl)\{g(Y,Z)X-g(X,Z)Y\}+\frac{1}{2}(3trl-r)\{g(Y,Z)\eta(X)\xi\\
	\label{2.5}&-g(X,Z)\eta(Y)\xi+\eta(Y)\eta(Z)X-\eta(X)\eta(Z)Y\},\\
	\label{2.6}QX=&\frac{1}{2}(r-trl)X+\frac{1}{2}(3trl-r)\eta(X)\xi.
\end{align}
Equation \eqref{2.6} is equivalent to
\begin{align}
	\label{2.7}Ric(X,Y)=\frac{1}{2}(r-trl)g(X,Y)+\frac{1}{2}(3trl-r)\eta(X)\eta(Y).
\end{align}
\section{Proof of Main Results}
\subsection{Proof of Theorem \ref{t1.1}}
\begin{proof}
	We distinguish the cases $(i)\,\,m=\infty$ and $(ii)\,\, 0<m<\infty$. The first case together with \eqref{1.2} entails that the associated Riemannian metric $g$ is Yamabe soliton. Hence the conclusion follows from the result of Venkatesha and Naik (see Theorem 2 in \cite{Venkatesha}). Next, we proceed to the case $(ii)$. The soliton vector field $V$ is collinear with the Reeb vector field $\xi$ implies $V=\sigma\xi$, for some smooth function $\sigma$. Taking covariant differentiation of this along $X$ and making use of \eqref{2.3} we ultimately have
	\begin{align}
		\label{3.1}\nabla_X V=(X \sigma)\xi-\sigma(\varphi X+\varphi hX).
	\end{align} 
	Observe that the equation \eqref{1.2} implies
	\begin{align}
		\label{3.2}\frac{1}{2}\{g(\nabla_XV,Y)+g(X,\nabla_YV)\}=(r-\lambda)g(X,Y)+\frac{\sigma^2}{m}\eta(X)\eta(Y).
	\end{align}
	In view of \eqref{3.1}, the equation \eqref{3.2} transform into
	\begin{align}
		\label{3.3}\frac{1}{2}\{(X\sigma)\eta(Y)+(Y\sigma)\eta(X)\}-\sigma g(\varphi hX,Y)=(r-\lambda)g(X,Y)+\frac{\sigma^2}{m}\eta(X)\eta(Y).
	\end{align}
	Substituting $X$ and $Y$ by $\xi$ in \eqref{3.3}, we obtain
	\begin{align}
		\label{3.4}(\xi \sigma)=(r-\lambda)+\frac{\sigma^2}{m}.
	\end{align}
	On the other hand, replacing $Y$ by $\xi$ in \eqref{3.3} and making use of the previous equation yields
	\begin{align}
		\label{3.5}D\sigma-(\xi \sigma)\xi=0.
	\end{align}
	Differentiating this and making use of \eqref{2.3} we obtain
	\begin{align*}
		\nabla_X D\sigma=X(\xi \sigma)\xi-(\xi \sigma)\{\varphi X+\varphi hX\}.
	\end{align*}
	Applying the Poincare lemma ($d^2=0$), we see that
 \[
 X(\xi \sigma)\eta(Y)-Y(\xi \sigma)\eta(X)+2(\xi \sigma)d\eta(X,Y)=0.
 \]
 Choosing $X,Y\perp\xi$ and since $d\eta$ is non-vanishing on any contact metric manifold, we have $(\xi \sigma)=0$. From this, the equation \eqref{3.5} shows that $\sigma$ is constant. Thus the equation \eqref{3.4} entails that the scalar curvature $r$ is also constant. Moreover, the equation \eqref{3.3} reduces to
	\begin{align*}
		-\sigma g(\varphi hX,Y)=(r-\lambda)g(X,Y)+\frac{\sigma^2}{m}\eta(X)\eta(Y).
	\end{align*}
	Differentiating this along $Z$ and making use of \eqref{2.3} we obtain
	\begin{align}
		-\sigma g((\nabla_Z\varphi h)X,Y)
		= (Zr)g(X,Y)-\frac{\sigma ^2}{m} & \left\{ \eta(X)g(Y,\varphi Z+\varphi hZ) \right. \nonumber \\
		&+ \left. \eta(Y)g(X,\varphi Z+\varphi hZ) \right\}.
	\end{align}
	Contracting this over $Y$ and $Z$ and since $r$ is constant, $tr\varphi=tr\varphi h=0$, we obtain
	\begin{align}
		\label{3.7}\sigma(div(\varphi h))X=0.
	\end{align}
	Remember that for any contact metric manifold $div(\varphi h)X=2n\eta(X)-g(Q\xi,X)$ (see\cite{Blair}). In view of $Q\varphi=\varphi Q$, equation~\eqref{2.4} and $\varphi \xi=0$, we have that $Q\xi=(trl)\xi$.  This yields $div(\varphi h)X=(2-trl)\eta(X)$. By virtue of this, the equation \eqref{3.7} entails that 
	\begin{align*}
		\sigma(2-trl)=0.
	\end{align*}
	Thus we have either $\sigma=0$, or $trl=2$. The first case shows that soliton vector field $V=0$, and this leads to a contradiction as $V$ is non-zero. So, we must have $trl=2$. By virtue of this, the equation \eqref{2.4} yields $trh^2=0$ and hence, since $h$ is symmetric, $h=0$ and so $M$ is Sasakian. Moreover, since $\xi$ is Killing and $\sigma$ is constant, the vector field $V(=\sigma\xi)$ is also Killing. This completes the proof of the Theorem. 
\end{proof} 
Before enter the proof of Theorem \ref{t1.2}, we give some key lemmas, which will be helpful to prove Theorem \ref{t1.2}.
\begin{lemma}
	Let $M$ be a 3-dimensional contact metric manifold with $Q\varphi=\varphi Q$. If $M$ admits a quasi Yamabe gradient soliton, then the following formulae are valid
	\begin{align}
		\label{3.8}(i)\,\, R(X,Y)Df&=(Xr)Y-(Yr)X+\frac{r-\lambda}{m}\{(Yf)X-(Xf)Y\}.\\
		\label{3.9}(ii)\,\, Ric(Y,Df)&=-2(Yr)+\frac{2(r-\lambda)}{m}(Yf).
	\end{align}
\end{lemma}
\begin{proof}
	We note that the equation \eqref{1.3} can be expressed as
	\begin{align}
		\label{3.10}\nabla_Y Df=(r-\lambda)Y+\frac{1}{m}g(Y,Df)Df,
	\end{align}
	where $D$ is the gradient operator of $g$. Differentiating \eqref{3.10} along an arbitrary vector field $X$, we obtain
	\begin{align*}
		\nonumber\nabla_X\nabla_Y Df=(Xr)Y+(r-\lambda)\nabla_X Y+\frac{1}{m}\{g(\nabla_X Y,Df)Df\\
		+g(Y,\nabla_X Df)Df+g(Y,Df)\nabla_X Df\}.
	\end{align*}
	Making use of previous equation and \eqref{3.10} in the well known expression of curvature tensor,
	\begin{align*}
		R(X,Y)=[\nabla_X,\nabla_Y]-\nabla_{[X,Y]},
	\end{align*}
	yields $(i)$. Next, contracting \eqref{3.8} over $X$ we get $(ii)$. This finishes the proof.
\end{proof}
\begin{lemma}\label{l3.2}\cite{BKS}
	Let $M$ be a contact metric manifold with a contact metric structure $(\varphi,\xi,\eta,g)$ such that $Q\varphi=\varphi Q$. Then the function $trl$ is constant everywhere on $M$ and $(\xi r)=0$. Further, if $trl=0$ then $M$ is flat.
\end{lemma}
\subsection{Proof of Theorem \ref{t1.2}}
\begin{proof}
	Taking inner product of \eqref{3.8} with $\xi$ we obtain
	\begin{align}
		\label{3.11}g(R(X,Y)Df,\xi)=(Xr)\eta(Y)-(Yr)\eta(X)+\frac{r-\lambda}{m}\{(Yf)\eta(X)-(Xf)\eta(Y)\}.
	\end{align}
	Substituting $Z$ by $\xi$ in \eqref{2.5} and then taking inner product of resultant equation with $\xi$ one gets
	\begin{align}
		\label{3.12}g(R(X,Y)\xi,Df)=(\frac{trl}{2})\{(Xf)\eta(Y)-(Yf)\eta(X)\}.
	\end{align}
	Comparing \eqref{3.12} with \eqref{3.11} we infer
	\begin{align*}
		(Xr)\eta(Y)-(Yr)\eta(X)+(\frac{r-\lambda}{m}-(\frac{trl}{2}))\{(Yf)\eta(X)-(Xf)\eta(Y)\}=0.
	\end{align*}
	It was proved in \cite{Huang} that the scalar curvature of any compact quasi Yamabe gradient soliton should be constant. By virtue of this, the previous equation reduces to
	\begin{align*}
		\frac{2(r-\lambda)-m(trl)}{2m}\{(Yf)\eta(X)-(Xf)\eta(Y)\}=0,
	\end{align*}
	which is equivalent to
	\begin{align}
		\label{3.13}\frac{2(r-\lambda)-m(trl)}{2m}\{Df-(\xi f)\xi\}=0.
	\end{align}
	Since $r$, $trl$ and $\lambda$ are constant, we have either $trl=\frac{2(r-\lambda)}{m}$, or $trl\neq\frac{2(r-\lambda)}{m}$. We now discuss the two cases separately.\\
	$\textbf{Case(i)}$. In this case, we have $trl=\frac{2(r-\lambda)}{m}$. Then from \eqref{3.10} it follows that
	\begin{align*}
		g(\nabla_\xi Df,\xi)=(r-\lambda)+\frac{1}{m}(\xi f)^2.
	\end{align*}
	We know that $g(\xi,Df)=(\xi f)$. Differentiation of this along $\xi$ and noting that $\nabla_\xi \xi=0$ (a consequence of \eqref{2.3}), we obtain $g(\nabla_\xi Df,\xi)=\xi(\xi f)$. From the above equation, one can find
	\begin{align}
		\label{3.14}\xi(\xi f)=(r-\lambda)+\frac{1}{m}(\xi f)^2.
	\end{align}
	As a result of \eqref{2.7} and \eqref{3.9}, one can get
	\begin{align*}
		(\frac{4(r-\lambda)}{m}-trl)(\xi f)=0.
	\end{align*}
	By virtue of $trl=\frac{2(r-\lambda)}{m}$, it follows from the foregoing equation that $(\xi f)=0$. Making use of this in \eqref{3.14} one can obtain $(r-\lambda)=0$. Utilization of this in $trl=\frac{2(r-\lambda)}{m}$ we find $trl=0$. From Lemma \ref{l3.2} we conclude that $M$ is flat.\\
	$\textbf{Case(ii)}$. When $trl\neq\frac{2(r-\lambda)}{m}$, from \eqref{3.13} we have $Df=(\xi f)\xi$. Differentiating this along $X$ and making use of \eqref{2.3} one can get $\nabla_X Df=X(\xi f)\xi-(\xi f)\{\varphi X+\varphi hX\}$. Applying Poincare lemma: $d^2=0$ gives
	\begin{align*}
		X(\xi f)\eta(Y)-Y(\xi f)\eta(X)+2(\xi f)d\eta(X,Y)=0.
	\end{align*} 
	Choosing $X,Y\perp \xi$ and we know that $d\eta$ is non-vanishing on any contact metric manifold, we have $(\xi f)=0$. From this we have $Df=0$. This
	shows that $f$ is constant, which in turn implies that $M$ admits a trivial quasi yamabe gradient soliton. This establishes the proof of the theorem.	
\end{proof}

\begin{remark}
	It is easy to seen that, the conclusion in Theorem~\ref{t1.2} still holds if we replace \textit{compactness} by \textit{constant scalar}.
\end{remark}

 \section*{Acknowledgments}
The authors are thankful to Department of Science and Technology, New Delhi for financial assistance to the Department of Mathematics, Kuvempu University under the FIST program (Ref. No. SR/FST/MS-I/2018-23(C)). Also, the authors wishes to thank the Referee for his/her constructive suggestions in improving the paper.

%%% References

\EditInfo{%
    June 24, 2019}{%
    July 03, 2020}{%
    Haizhong Li}

\end{paper}